\documentclass[10pt]{article}
\usepackage{amssymb, amsmath}
\usepackage{float,epsfig}

\begin{document}

\newtheorem{theorem}{\bf Theorem}[section]
\newtheorem{proposition}[theorem]{\bf Proposition}
\newtheorem{definition}[theorem]{\bf Definition}
\newtheorem{corollary}[theorem]{\bf Corollary}
\newtheorem{example}[theorem]{\bf Example}
\newtheorem{exam}[theorem]{\bf Example}
\newtheorem{remark}[theorem]{\bf Remark}
\newtheorem{lemma}[theorem]{\bf Lemma}
\newcommand{\nrm}[1]{|\!|\!| {#1} |\!|\!|}

\newcommand{\ba}{\begin{array}}
\newcommand{\ea}{\end{array}}
\newcommand{\von}{\vskip 1ex}
\newcommand{\vone}{\vskip 2ex}
\newcommand{\vtwo}{\vskip 4ex}
\newcommand{\dm}[1]{ {\displaystyle{#1} } }

\newcommand{\be}{\begin{equation}}
\newcommand{\ee}{\end{equation}}
\newcommand{\beano}{\begin{eqnarray*}}
\newcommand{\eeano}{\end{eqnarray*}}
\newcommand{\inp}[2]{\langle {#1} ,\,{#2} \rangle}
\def\bmatrix#1{\left[ \begin{matrix} #1 \end{matrix} \right]}
\def \noin{\noindent}
\newcommand{\evenindex}{\Pi_e}



\def \R{{\mathbb R}}
\def \C{{\mathbb C}}
\def \K{{\mathbb K}}
\def \J{{\mathbb J}}
\def \Lb{\mathrm{L}}

\def \T{{\mathbb T}}

\def \norm{\nrm{\cdot}\equiv \nrm{\cdot}}

\def \H{{\mathbb H}}
\def \L{{\mathbb L}}

\def \S{\mathbb{S}}
\def \sigmin{\sigma_{\min}}
\def \elam{\Lambda_{\epsilon}}

\def \N{{\mathbb N}}
\def \tr{\mathrm{Tr}}
\def \lam{\lambda}
\def \sig{\sigma}
\def \sign{\mathrm{sign}}
\def \ep{\epsilon}
\def \diag{\mathrm{diag}}

\def \rank{\mathrm{rank}}
\def \pf{{\bf Proof: }}
\def \dist{\mathrm{dist}}
\def \rar{\rightarrow}

\def \pf{{\bf Proof: }}
\def \Re{\mathsf{Re}}
\def \Im{\mathsf{Im}}
\def \re{\mathsf{re}}
\def \im{\mathsf{im}}

\def \sym{\mathsf{sym}}
\def \sksym{\mathsf{skew\mbox{-}sym}}
\def \odd{\mathrm{odd}}
\def \even{\mathrm{even}}
\def \herm{\mathsf{Herm}}
\def \skherm{\mathsf{skew\mbox{-}Herm}}
\def \str{\mathrm{ Struct}}
\def \eproof{$\blacksquare$}
\def \proof{\noin\pf}

\def \cA{\mathcal{A}}
\def \cG{\mathcal{G}}
\def \F{{\mathcal F}}
\def \cH{\mathcal{H}}
\def \cJ{\mathcal{J}}
\def \tr{\mathrm{Tr}}
\def \range{\mathrm{Range}}
\def \adj{\star}
\def \argmin{\arg\!\min}


\title{\sc  Structured inverse least-squares problem for structured matrices}
\author{ Bibhas Adhikari\thanks{Centre for System Science,
IIT Jodhpur, India, E-mail:
bibhas@iitj.ac.in } \, and Rafikul Alam\thanks{Corresponding author: Department of Mathematics,
IIT Guwahati, India, E-mail: rafik@iitg.ernet.in, Fax:~+91-361-2690762/2582649. }  }
\date{}

\maketitle
\thispagestyle{empty}

{\small \noin{\bf Abstract.}  Given a pair of matrices $X$ and $B$ and an appropriate class of  structured matrices $\S,$  we provide a complete solution of the structured inverse least-squares problem $\min\limits_{A \in \S}\|AX-B\|_F.$ Indeed, we determine all solutions of the structured inverse least squares problem as well as those solutions which have the smallest norm. We show that there are infinitely many smallest norm solutions of the least squares problem for the spectral norm whereas  the smallest norm solution is unique for the Frobenius norm. }

\vone \noin{\bf Keywords.} Structured matrices, structured backward errors, Jordan and Lie algebras,  structured inverse least-squares.

\vone\noin {\bf AMS subject classification(2000):}  15A24,  65F15, 65F18,  15A18.

\section{Introduction}  Let $\S$ be a class of structured matrices in $\K^{n\times n},$ where $ \K= \R$ or $\C.$ We consider  $\S$ to be either a Jordan algebra or a Lie algebra associated with an appropriate scalar product on $\K^n.$ This provides a general setting that encompasses important classes of structured matrices such as Hamiltonian, skew-Hamiltonian, symmetric, skew-symmetric, pseudosymmetric, persymmetric, Hermitian, Skew-Hermitian, pseudo-Hermitian, pseudo-skew-Hermitian, to name only a few, see~\cite{fran6}.

Given  a pair of matrices $X$ and $B$ in $\K^{n\times p}$ and a class $\S$ of structured matrices,  finding a matrix $A \in \S$ such that $AX= B$  is known as the {\em structured mapping problem}~\cite{AlaA, fran6, khatri}.
%
Structured mapping problem  arises in many applications~\cite{anto,griv, AlaBKMM, tiss03, fran6}. For example,  the passivation problem for a stable linear time-invariant (LTI) control  system
\be \label{paslti}
\begin{array}{lcl}
  \dot x & = & \mathcal{A} x + \mathcal{B} u,\ x(0)=0, \\
  y & = & \mathcal{C} x + \mathcal{D} u,
\end{array}
\ee
gives rise to a structured mapping problem for Hamiltonian matrices~\cite{AlaBKMM}, where
$\mathcal{A}\in \mathbb K^{n\times n}$,
$\mathcal{B}\in \mathbb K^{n\times p}$, $\mathcal{C}\in \mathbb K^{p\times n}$,
$\mathcal{D}\in \mathbb K^{p\times p},$   $u$ is the input, $x$ is the state and $y$ is the output. 
A complete solution of the structured mapping problem has been provided in~\cite{AlaA} when the class of structured matrices $\S$ is either a Jordan algebra or a Lie algebra associated with an orthosymmetric scalar product on $\K^n,$ see also~\cite{fran6, tiss03}. It is shown that the matrix equation $AX =B$ admits a structured solution $A \in \S$ if and only if  $X$ and $B$ satisfy appropriate conditions~\cite{AlaA, Ba:thesis}.



We consider a more general problem,  namely, the {\em structured inverse least-squares problem} (SILSP) which minimizes  $\|AX-B\|_F$ when  $ A$ varies  in $\S,$ where $\|\cdot\|_F$ is the Frobenius norm. Note that if $AX = B$ admits a solution $A \in \S$ then $A$ is also a solution of the SILSP.

%

\vone
\noin{\sc Problem-I.} {\bf (Structured inverse least-squares problem)} Let $\S \subset \K^{n\times n}$ be a class of structured matrices and let $ X, B \in \K^{n\times p}.$ Set \beano  \rho^{\S}(X,B) &:=& \inf\{ \|AX-B\|_F: A \in \S\}, \\
  \sigma^{\S}(X, B) &:= & \inf\{ \|A\| : A \in  \S \mbox{ and } \|AX-B\|_F = \rho^{\S}(X,B) \}. \eeano  Determine all matrices in $\S$ such that $\|AX-B\|_F = \rho^{\S}(X,B).$ Also determine all optimal solutions $A_o$  in $\S$ such that $\|A_oX-B\|_F = \rho^{\S}(X,B)$ and $\|A_o\| =\sigma^{\S}(X, B),$ where $\|\cdot\|$ is either the spectral norm or the Frobenius norm.

\vone
Structured inverse least-squares problem aries in many applications such as in particle
physics and geology, inverse problems of vibration theory, inverse Sturm-Liouville problem, control theory, and in multidimensional
approximation, see~\cite{dai,fau, wood} and references therein. Structured inverse least-squares problem has been considered for positive semidefinite solutions in~\cite{wood} and for skew-Hermitian Hamiltonian solutions in ~\cite{zha:liu}.

 We provide a complete solution of the structured inverse least-squares problem  when the class $\S \subset \K^{n\times n}$ of structured matrices is either a Jordan algebra or a Lie algebra associated with an orthosymmetric scalar product on $\K^n.$  We characterize solutions of the structured inverse least-squares problem and show that for the spectral norm there are infinitely many optimal structured solutions (that is, structured solutions having the smallest spectral norm). On the other hand, we show that for the Frobenius norm the optimal structured solution (that is, structured solution having the smallest Frobenius norm) is unique.  We determine all optimal solutions of the structured inverse least-squares problem ({\sc Problem-I}) for the spectral  and the  Frobenius norms.

 \vone

\noindent{\bf Notation.} Let $\K^{m\times n}$ denote the set of all $m$-by-$n$ matrices with entries in $\K,$ where $\K = \R$ or $\C.$
We denote the transpose of a matrix $X \in \K^{m\times n}$ by $X^T,$  the conjugate transpose by $X^H$ and the complex conjugate by $\bar X.$  We often write $A^*$ for $* \in \{T, H\}$ to denote the transpose $A^T$ or the conjugate transpose $A^H.$ We denote the spectral norm and the Frobenius norm of $X$ by  $\|X\|_2$ and $ \|X\|_F,$ respectively, and are given by $$ \|X\|_2 := \max\{\|Xu\|_2 : \|u\|_2 =1\}  \, \mbox{ and } \,  \|X\|_F := \sqrt{\tr(X^HX)}.$$


\section{Structured matrices}
We now briefly define structured matrices that we consider in this paper; see~\cite{fran6} for further details.
Let  $ M \in \K^{n\times n}$ be unitary. Assume further that $M$ is either symmetric  or skew-symmetric or Hermitian or skew-Hermitian.
Define the scalar product $\inp{\cdot}{\cdot}_M : \K^n\times \K^n \rar \K$ by \be \label{sprod} \inp{x}{y}_M := \left\{\ba{cc} y^TMx, & \mbox{ bilinear form, }\\
y^HMx, & \mbox{ sesquilinear form.} \ea\right.\ee Then for $A \in \K^{n\times n}$ there is a unique adjoint operator $A^\star$ relative to the scalar product (\ref{sprod}) such that
 $\inp{Ax}{y}_M= \inp{x}{A^\star y}_M $ for all $x, y \in \K^n.$ The adjoint $A^\star$ is explicitly given by
\be \label{adj} A^{\star} = \left\{ \ba{cc} M^{-1}A^{T}M, & \mbox{bilinear form,}\\
               M^{-1}A^HM, & \mbox{ sesquilinear form.}
               \ea\right.\ee
Consider the Lie algebra $\L$ and the Jordan algebra $\J$ associated with the scalar product (\ref{sprod}) given by
\be \label{alg} \L :=\{ A\in \K^{n \times n}: A^{\star}=-A\} \mbox{ and } \J :=\{ A\in \K^{n \times n}: A^{\star}=A\}.\ee
In this paper, we consider $\S = \L $ or $\S= \J,$ and refer to the matrices in $\S$ as {\em structured matrices.}
The Jordan and Lie algebras  so defined provide a general framework for analyzing many important classes of structured matrices including Hamiltonian, skew-Hamiltonian, symmetric, skew-symmetric, pseudosymmetric, persymmetric, Hermitian, Skew-Hermitian, pseudo-Hermitian, pseudo-skew-Hermitian matrices, to name only a few, see~(Table~2.1, \cite{fran6}).

For the rest of the paper, we set \beano \sym &:=&
\{A\in \K^{n \times n} : A^T = A\}, \, \,\, \sksym  \,\, := \,\, \{A\in \K^{n \times n} : A^T = -A\},  \\ \herm &:=& \{A\in \C^{n \times n} :
A^H = A\}, \, \,\, \skherm \,\, :=\,\,  \{A\in\C^{n \times n } : A^H = -A\}. \eeano
Also define the set $M\S$ by  $ M\S :=\{ M A : A \in \S\}.$ Then, in view of (\ref{adj}) and (\ref{alg}), it follows that
\be \label{str4} \S \in \{ \J, \L\} \Longleftrightarrow M \S \in \{ \sym, \, \sksym, \, \herm, \, \skherm\}.\ee
This shows that the four classes of structured matrices, namely, symmetric, skew-symmetric, Hermitian and skew-Hermitian matrices are prototypes of more general structured matrices belonging to the Jordan and Lie algebras given in (\ref{alg}).

Let $A, B, C$ and $D$ be matrices.  Then the matrix $T :=\bmatrix{ A & C\\
B & D}$ is called a {\bf dilation} of $A.$ The norm preserving
dilation problem is then stated as follows. Given matrices $A, B, C$
and a positive number \be  \label{mu} \mu \geq \max\left( \left\|\bmatrix{ A\\
B}\right\|_2, \, \left\|\bmatrix{ A & C}\right\|_2\right),\ee find all possible $D$ such that $\left\|\bmatrix{ A &
C\\ B & D}\right\|_2 \leq \mu.$ 

\begin{theorem}[Davis-Kahan-Weinberger, \cite{D}] \label{dkw}
Let $A, B, C$ be given matrices. Then for any positive number
$\mu$ satisfying (1),
  there exists $D$ such that $\left\|\bmatrix{  A & C \\
  B & D
}\right\|_2 \leq \mu.$ Indeed, those $D$ which have this property
are exactly those of the form $$D= - KA^HL +
\mu(I-KK^H)^{1/2}Z(I-L^HL)^{1/2},$$ where $K^H :=(\mu^2 I -A^HA
)^{-1/2}B^H,~ L :=(\mu^2 I -AA^H)^{-1/2}C$ and $Z$ is an
arbitrary contraction, that is, $\|Z\|_2 \leq 1.$
\end{theorem}

We mention that when $(\mu^2 I-A^HA)$ is singular, the
inverses in $K^H$ and $L$ are replaced by their Moore-Penrose
pseudo-inverses~(see, \cite{meinguet}). An interesting fact about Theorem~\ref{dkw} is that if $T(D)$ is symmetric (resp., skew-symmetric, Hermitian, skew-Hermitian) then the solution matrices $D$ is symmetric (resp., skew-symmetric, Hermitian,  skew-Hermitian).

Let $\S \in \{\sym, \sksym, \herm, \skherm \}$ and  $A \in \S$  be given by  \be\label{fronorm} A := \bmatrix{A_{11} & \pm A_{12}^* \\
A_{12} & A_{22}}. \mbox{  Then  } \|A\|_F = \left(2\left\|
\bmatrix{A_{11} \\ A_{12}} \right\|_F^2 - \|A_{11}\|_F^2 +
\|A_{22}\|_F^2\right)^{1/2}\ee where $ * \in \{T, H\}.$  We repeatedly use this fact in the sequel.

\section{Structured inverse least-squares problem }
Consider the set of structured matrices $\S \in \{\J, \L\},$  where $\J$ is the Jordan algebra and  $\L$ is the Lie algebra  given in (\ref{alg}).
Let $ X$ and $B$ be a pair of matrices in $\K^{n\times p}.$ Then, since  $M$ is unitary and $\|A X - B\|_F = \|MA X - MB\|_F $ for $A \in \S,$ it follows that  $$ \min_{A \in \S} \|AX - B\|_F = \min_{ G \in M\S}\|GX-MB\|_F.$$   Consequently, in view of (\ref{str4}), we only need to consider {\sc Problem-I} for the special case when  $\S \in \{ \sym, \, \sksym, \, \herm, \, \skherm\}.$ Further, if $\cG(X, B) $ denotes a solution of {\sc Problem-I} for the structured matrices $M\S$   then  $A := M^{-1} \cG(X, MB)$ provides a solution of {\sc Problem-I} for the general case when $ \S \in \{ \J, \L\}.$ In other words, we have \be \label{solmap}  \cG(X, B)  \in  \argmin\limits_{A \in M\S} \|AX -B\|_F  \Longleftrightarrow M^{-1}\cG(X, MB) \in \argmin\limits_{A \in \S} \|AX -B\|_F.\ee Furthermore, if $\cG(X, B)$ is an optimal solution (has the smallest norm) then so is $M^{-1}\cG(X, MB)$ and that \be \label{eqlrho} \rho^{\S}(X, B) =\rho^{M\S}(X, MB).\ee
Note that $\cG(X, B)$ is skew-Hermitian solution of $\min_A \|AX-B\|_F$  if and only of $ -i \cG(X, \, iB)$ is a Hermitian solution of $\min_A\|AX-B\|_F.$ Consequently, we only need to solve {\sc Problem-I} for the class of structured matrices $\S \in \{ \sym, \, \sksym, \herm\}.$

\vone


Let $(X,B)\in \K^{n\times k}\times \K^{n\times k}.$ Suppose that  $\rank(X)=r$ and consider the``trimmed"  SVD $X=U_1\Sigma_1 V_1^H,$ where $\Sigma_1=\diag\{\sigma_1, \hdots, \sigma_r\}.$ Let $D \in \C^{r\times r}$ be given by  $D_{ij}:=1/(\sigma_i^2 + \sigma_j^2).$ For the rest of this section, define
{\small
\begin{eqnarray*}   \mathcal{F}_\pm^T(X,B) &:=& \overline{U}_1 [D \circ (\Sigma_1 U_1^TBV_1 \pm V_1^TB^TU_1\Sigma_1)]U_1^H
 \pm (BX^\dagger)^T(I - XX^\dagger) + (I-XX^\dagger)^TBX^\dagger \\ \mathcal{F}_+^H(X,B) &:=& U_1 [D \circ (\Sigma_1 U_1^HBV_1 + V_1^HB^HU_1\Sigma_1)]U_1^H
 + (BX^\dagger)^H(I - XX^\dagger) + (I-XX^\dagger)^HBX^\dagger \end{eqnarray*} } where $\circ$ denotes the Hadamard product.

We repeatedly employ the following elementary result in the sequel.
\begin{lemma}\label{lemma1}
Let $\alpha$ and $ \beta$ be real numbers and $ b_1, b_2 \in \C.$ Then $\min_{x \in
\C} (|x\alpha -b_1|^2 + |x\beta - b_2|^2)$ is attained at $x =
\dfrac{\alpha b_1 + \beta b_2}{\alpha^2 + \beta^2}.$
\end{lemma}

The next result determines $\rho^{\S}(X, B)$ and characterizes $ \argmin\limits_{A \in \S} \|AX -B\|_F$.

\begin{theorem}\label{silsp:sym}
Let $(X,B)\in \K^{n\times p}\times \K^{n\times p}$ and $\S\in \{\sym, \sksym, \herm\}.$ Suppose that $\rank(X)=r$ and consider the SVD $X=U\Sigma V^H.$
Partition $U=[U_1 \,\, U_2]$ and $ V=[V_1 \,\, V_2],$ where $U_1,V_1\in\K^{n\times r}.$ Let $\Sigma_1:=\Sigma(1:r,1:r)=\diag(\sigma_1, \hdots, \sigma_r).$ Then  we have $$\rho^{\S}(X, B) =\left\{%
\begin{array}{ll}
\sqrt{ \|(D \circ (\Sigma_1 U_1^*BV_1 +
V_1^*B^*U_1 \Sigma_1))\Sigma_1 - U_1^*BV_1\|_F^2
 + \|BV_2\|_F^2}, \\  \hfill{\mbox{if} \,\, \S\in\{\sym,\herm\}}, \\ & \\
\sqrt{ \| (D
\circ (\Sigma_1 U_1^TBV_1 - V_1^TB^TU_1 \Sigma_1))\Sigma_1 -
U_1^TBV_1\|_F^2
 + \|BV_2\|_F^2 } , \\ \hfill{\mbox{if} \,\, \S=\sksym}, \\
\end{array}%
\right. $$ where $*=T$ when $\S=\sym,$ and $*=H$ when $\S=\herm.$ Here $ D \in \C^{r\times r}$ is given by $ D_{ij} := \dfrac{1}{\sigma_i^2 +
 \sigma_j^2}$ and $\circ$ denotes the Hadamard product.

Further,  $A_o\in \S$ and  $\|A_oX-B\|_F=\rho^{\S}(X, B)$ if and only if
$A_o$ is of the form  $$A_o =\left\{%
\begin{array}{ll}
\mathcal{F}_+^T(X,B) + (I-XX^\dagger)^T Z (I-XX^\dagger), \,\, \mbox{if} \,\, \S=\sym \\
\mathcal{F}_-^T(X,B) + (I-XX^\dagger)^T Z (I-XX^\dagger),\,\, \mbox{if} \,\, \S=\sksym \\
\mathcal{F}_+^H(X,B) + (I-XX^\dagger)^H Z (I-XX^\dagger), \,\, \mbox{if} \,\, \S=\herm \\
\end{array}%
\right. $$ for some $Z\in\S$.


\end{theorem}

\noin\pf Suppose that $\S=\sym.$ Then we have  $$ A = \overline{U} \, \overline{U}^HAUU^H= \bmatrix{\overline{U}_1 & \overline{U}_2} \bmatrix{A_{11} & A_{12}^T \\
A_{12} & A_{22}} \bmatrix{U_1^H \\ U_2^H} $$ where $A_{11} = A^T_{11} \in \C^{r\times r}$  and, $
A_{12}$ and $ A_{22} = A_{22}^T$ are matrices of compatible sizes. Then
\beano \|AX-B\|_F^2 &=& \|U^TAUU^HX - U^TB\|_F^2  =
\left\|\bmatrix{A_{11} & A_{12}^T \\ A_{12} & A_{22}} \bmatrix{
U_1^HX
\\ 0} - \bmatrix{U_1^TB \\ U_2^TB}\right\|_F^2 \\ &=& \|A_{11}U_1^HX - U_1^TB\|_F^2 +
 \| A_{12}U_1^HX - U_2^TB\|_F^2. \eeano
Note that   \beano  \| A_{12}U_1^HX - U_2^TB\|_F^2
 &=&  \|A_{12}U_1^HU\Sigma V^H - U_2^TB\|_F^2
 =  \|A_{12}U_1^HU\Sigma - U_2^TBV \|_F^2 \\
 &=&  \|A_{12}\Sigma_1 - U_2^TBV_1\|_F^2 + \|U_2^TBV_2\|_F^2  \geq \|U_2^TBV_2\|_F^2 \eeano
and the minimum is attained when  $A_{12}\Sigma_1 - U_2^TBV_1=0,$ that is, when $A_{12}
= U_2^TBV_1\Sigma_1^{-1}.$
Also  note that
$$ \|A_{11}U_1^HX - U_1^TB\|_F^2 =
  \|A_{11}\Sigma_1 - U_1^TBV_1\|_F^2 + \|U_1^TBV_2\|_F^2 $$
is minimized when  $\|A_{11}\Sigma_1 - U_1^TBV_1\|_F^2$ is minimized over the symmetric matrices $A_{11}.$


Suppose that $A_{11} = [a_{ij}] $ and $ U_1^TBV_1 = [b_{ij}]. $ Then by Lemma~\ref{lemma1}
 $$\|A_{11}\Sigma_1 - U_1^TBV_1\|_F^2
= \sum_{i < j,i=1}^r \sum_{j=1}^r (|a_{ij}\sigma_j -
b_{ij}|^2 + |a_{ij}\sigma_i - b_{ji}|^2)$$
is minimized when $$a_{ij} = \frac{\sigma_j b_{ij} + b_{ji}\sigma_i}{\sigma_i^2 +
\sigma_j^2}, \, a_{ij} = a_{ji}, \mbox{ for } i,j =1:r.$$ This yields
$$A_{11} = D \circ (U_1^TBV_1\Sigma_1 + \Sigma_1V_1^TB^TU_1).$$
Hence  \beano \rho^{\S}(X,B) &=& \sqrt{ \| (D \circ (\Sigma_1 U_1^TBV_1 +
V_1^TB^TU_1 \Sigma_1)) \Sigma_1 - U_1^TBV_1\|_F^2
 + \|U_1^TBV_2\|_F^2 + \|U_2^T B
V_2\|_F^2} \\ &=& \sqrt{ \| (D \circ (\Sigma_1 U_1^TBV_1 +
V_1^TB^TU_1 \Sigma_1)) \Sigma_1 - U_1^TBV_1\|_F^2
 + \|BV_2\|_F^2}\eeano as desired.

Now, substituting $A_{11}$ and $A_{12},$ we have   \be\label{silsp:symaa} A = \overline{U}
 \bmatrix{D\circ (\Sigma_1 U_1^TBV_1 + V_1^TB^TU_1 \Sigma_1) &
 \Sigma_1^{-1}V_1^TB^TU_2 \\ U_2^TBV_1\Sigma_1^{-1} & A_{22}} U^H \ee which upon simplification  yields \beano A  &=& \overline{U}_1[D\circ
 (\Sigma_1 U_1^TBV_1 + V_1^TB^TU_1 \Sigma_1)]U_1^H + (BX^\dagger)^T(I - XX^\dagger) +
 (I - XX^\dagger)^TBX^\dagger  \\ &&+ (I - XX^\dagger)^T Z (I -
 XX^\dagger), \, \mbox{ for some } Z^T = Z \in \C^{n \times n}, \eeano as desired.

 Next, suppose that  $\S=\sksym.$ Then we have  $$ A = \bmatrix{\overline{U}_1 & \overline{U}_2} \bmatrix{A_{11} & -A_{12}^T \\
A_{12} & A_{22}} \bmatrix{U_1^H \\ U_2^H} $$ where $A_{11} = - A_{11}^T \in \C^{r\times r}$ and,
$A_{12}$ and $A_{22} = -A_{22}^T$ are matrices of compatible sizes.  Then as before
 $$\|AX-B\|_F^2 = \|A_{11}U_1^HX - U_1^TB\|_F^2 +
 \| A_{12}U_1^HX - U_2^TB\|_F^2$$ is minimized when  $A_{12}=U_2^TBV_1\Sigma_1^{-1}$ and $\|A_{11}U_1^HX - U_1^TB\|_F^2$ is minimized over all skew-symmetric matrices $A_{11}.$ Again  by Lemma~\ref{lemma1}, a minimizer is given by $$a_{ij}=\frac{\sigma_jb_{ij}-\sigma_ib_{ji}}{\sigma_i^2+\sigma_j^2}, \,\, a_{ji}=-a_{ij} \mbox{ for } i,j=1:r, $$ where $U_1^TBV_1=[b_{ij}],$ that is,
  $A_{11}=D \circ (U_1^TBV_1\Sigma_1 - \Sigma_1V_1^TB^TU_1). $ Consequently, we have
$$ \rho^{\S}(X, B) =  \sqrt{ \| (D \circ (\Sigma_1 U_1^TBV_1 -
V_1^TB^TU_1 \Sigma_1)) \Sigma_1 - U_1^TBV_1\|_F^2
 + \|BV_2\|_F^2}$$ and
 \begin{equation}\label{eqn2:sksym} A =\overline{U}
 \bmatrix{D\circ (\Sigma_1 U_1^TBV_1 - V_1^TB^TU_1 \Sigma_1) &
 -\Sigma_1^{-1}V_1^TB^TU_2 \\ U_2^TBV_1\Sigma_1^{-1} & A_{22}} U^H\end{equation} which upon simplification yields the desired form of $A.$

Finally, suppose that $\S=\herm.$ Then  we have $$A= UU^HAUU^H= [U_1 \,\, U_2]\bmatrix{A_{11} & A_{12}^H \\ A_{12} & A_{22}}\bmatrix{U_1^H \\ U_2^H},$$
where $A_{11}^H = A_{11} \in \C^{r\times r}$ and, $ A_{12}$ and $ A_{22}^H = A_{22}$ are matrices of compatible sizes. Again note that
 $$\|AX-B\|_F^2 = \|A_{11}U_1^HX - U_1^HB\|_F^2 +
 \| A_{12}U_1^HX - U_2^HB\|_F^2$$ is minimized when $A_{12}=U_2^HBV_1\Sigma_1^{-1}$ and $\|A_{11}U_1^HX - U_1^HB\|_F^2$ is minimized over all Hermitian matrices  $A_{11}.$
In view of Lemma~\ref{lemma1}, a minimizer is given by $$a_{ij}=\frac{\sigma_jb_{ij}+\sigma_i\overline{b}_{ji}}{\sigma_i^2+\sigma_j^2}, \,\, a_{ji}=\overline{a}_{ij} \mbox{ for }  i,j=1:r, $$ where $U_1^HBV_1=[b_{ij}].$ In other words,   $A_{11}=D \circ (U_1^HBV_1\Sigma_1 + \Sigma_1V_1^HB^HU_1).$
Thus, we have
$$\rho^{\S}(X, B) = \sqrt{ \| (D \circ (\Sigma_1 U_1^HBV_1 +
V_1^HB^HU_1 \Sigma_1)) \Sigma_1 - U_1^HBV_1\|_F^2
 + \|BV_2\|_F^2}$$
and  \begin{equation}\label{eqn2:herm} A =U
 \bmatrix{D\circ (\Sigma_1 U_1^HBV_1 + V_1^HB^HU_1 \Sigma_1) &
 \Sigma_1^{-1}V_1^HB^HU_2 \\ U_2^HBV_1\Sigma_1^{-1} & A_{22}} U^H\end{equation} which upon simplification yields the desired form of $A.$
This completes the proof. $\blacksquare$

\vone

The next result characterizes solutions in $ \argmin\limits_{A \in \S} \|AX -B\|_F$ which have the smallest norms and determines the norms.

\begin{theorem}\label{cor:silsp:sym}
Let $(X,B)\in \K^{n\times p}\times \K^{n\times p}$ and $\S\in \{\sym, \sksym, \herm\}.$ Suppose that $\rank(X)=r$ and consider the SVD $X=U\Sigma V^H.$
Partition $U=[U_1 \,\, U_2]$ and $ V=[V_1 \,\, V_2],$ where $U_1,V_1\in\K^{n\times r}.$ Let $\Sigma_1:=\Sigma(1:r,1:r)=\diag(\sigma_1, \hdots, \sigma_r).$

\begin{enumerate} \item {\bf Frobenius norm:} We have
$$\sigma^{\S}(X, B) =\left\{%
\begin{array}{ll}
\sqrt{ \|D \circ (\Sigma_1 U_1^*BV_1 +
V_1^*B^*U_1 \Sigma_1)\|_F^2
 + 2\|U_2^*BV_1\Sigma_1^{-1}\|_F^2}, \\  \hfill{\mbox{if} \,\, \S\in\{\sym,\herm\}}, \\ & \\
\sqrt{ \|D
\circ (\Sigma_1 U_1^TBV_1 - V_1^TB^TU_1 \Sigma_1)\|_F^2
 + 2\|U_2^TBV_1\Sigma_1^{-1}\|_F^2 } , \\ \hfill{\mbox{if} \,\, \S=\sksym}, \\
\end{array}%
\right. $$ where $*=T$ when $\S=\sym,$ and $*=H$ when $\S=\herm.$ Here $ D \in \C^{r\times r}$ is given by $ D_{ij} := \dfrac{1}{\sigma_i^2 +
 \sigma_j^2}$ and  $\circ$ denotes the Hadamard product.

There is a unique $A_o\in \S$ such that $\|A_oX-B\|_F=\rho^\S(X,B)$ and $\|A_o\|_F=\sigma^{\S}(X,B)$ and is given by $$A_o =\left\{%
\begin{array}{ll}
\mathcal{F}_+^T(X,B),  &  \mbox{ if } \,\, \S=\sym, \\
\mathcal{F}_-^T(X,B), &    \mbox{ if } \,\, \S=\sksym, \\
\mathcal{F}_+^H(X,B), &   \mbox{ if } \,\, \S=\herm. \\
\end{array}%
\right. $$

\item {\bf Spectral norm:} We have  $$\sigma^{\S}(X, B) =\left\{%
\begin{array}{ll}
\left\| \bmatrix{D \circ (\Sigma_1 U_1^*BV_1 + V_1^*B^*U_1 \Sigma_1) \\ U_2^*BV_1\Sigma_1^{-1}} \right\|_2 =:\mu_1,  & \mbox{ if }  \S\in\{\sym,\herm\}, \\ \\
\left\| \bmatrix{D \circ (\Sigma_1 U_1^TBV_1 - V_1^TB^TU_1 \Sigma_1) \\ U_2^TBV_1\Sigma_1^{-1}} \right\|_2 =:\mu_2,  & \mbox{ if} \,\, \S=\sksym, \\
\end{array}%
\right. $$ where $*=T$ when $\S=\sym,$ and $*=H$ when $\S=\herm.$
\end{enumerate}

Further,  $A_o \in \S$ such that $\|A_oX-B\|_F=\rho^\S(X,B)$ and $\|A_0\|_2=\sigma^{\S}(X,B)$ if and only if $A_o$ is of the form $$A_o =\left\{%
\begin{array}{ll}
\mathcal{F}_+^T(X,B) -(I-XX^\dagger)^TK [\overline{D \circ (\Sigma_1 U_1^*BV_1 + V_1^*B^*U_1 \Sigma_1)}]K^T(I-XX^\dagger)+\phi(Z), \\ \hfill{ \,\, \mbox{if} \,\, \S=\sym}, \\  & \\
\mathcal{F}_-^T(X,B) -(I-XX^\dagger)^TK [\overline{D \circ (\Sigma_1 U_1^*BV_1 - V_1^*B^*U_1 \Sigma_1)}]K^T(I-XX^\dagger)+\phi(Z), \\ \hfill{ \,\, \mbox{if} \,\, \S=\sksym}, \\ & \\
\mathcal{F}_+^H(X,B) -(I-XX^\dagger)^HK [D \circ (\Sigma_1 U_1^*BV_1 + V_1^*B^*U_1 \Sigma_1)]K^H(I-XX^\dagger)+\phi(Z), \\ \hfill{ \,\, \mbox{if} \,\, \S=\herm}, \\
\end{array}%
\right. $$ where $$K =\left\{%
\begin{array}{ll}
BX^\dagger U_1 [\mu_1^2I - (\overline{D \circ (\Sigma_1 U_1^TBV_1 + V_1^TB^TU_1 \Sigma_1)})(D \circ (\Sigma_1 U_1^TBV_1 + V_1^TB^TU_1 \Sigma_1))]^{-1/2}, \\ \hfill{ \,\, \mbox{if} \,\, \S=\sym}, \\ &\\
BX^\dagger U_1 [\mu_2^2I - (\overline{D \circ (\Sigma_1 U_1^*BV_1 - V_1^*B^*U_1 \Sigma_1)})(D \circ (\Sigma_1 U_1^*BV_1 - V_1^*B^*U_1 \Sigma_1))]^{-1/2}, \\ \hfill{ \,\, \mbox{if} \,\, \S=\sksym}, \\ & \\
BX^\dagger U_1 [\mu_1^2I - (D \circ (\Sigma_1 U_1^HBV_1 + V_1^HB^HU_1 \Sigma_1))(D \circ (\Sigma_1 U_1^HBV_1 + V_1^HB^HU_1 \Sigma_1))]^{-1/2}, \\
\hfill{ \,\, \mbox{if} \,\, \S=\herm}, \\
\end{array}%
\right. $$ and $$\phi(Z) =\left\{%
\begin{array}{ll}
 \mu_1\overline{U}_2 [I-U_2^TKK^H\overline{U}_2]^{1/2}Z [I-U_2^H\overline{K}K^TU_2]^{1/2}U_2^H, & Z=Z^T, \|Z\|_2\leq 1\\ & \hfill{ \,\, \mbox{if} \,\, \S=\sym,} \\ & \\
\mu_2\overline{U}_2 [I+U_2^TKK^H\overline{U}_2]^{1/2}Z [I+U_2^H\overline{K}K^TU_2]^{1/2}U_2^H, & Z=-Z^T, \|Z\|_2\leq 1\\ & \hfill{ \,\, \mbox{if} \,\, \S=\sksym,} \\ &\\
\mu_1U_2 [I-U_2^HKK^HU_2]^{1/2}Z [I-U_2^HKK^HU_2]^{1/2}U_2^H, & Z=Z^H, \|Z\|_2\leq 1\\ & \hfill{ \,\, \mbox{if} \,\, \S=\herm.} \\
\end{array}%
\right. $$
\end{theorem}

\noin\pf Suppose that $\S=\sym.$ Then by Theorem~\ref{silsp:sym} and (\ref{silsp:symaa}) we have  $$A_o = \overline{U}
 \bmatrix{D\circ (\Sigma_1 U_1^TBV_1 + V_1^TB^TU_1 \Sigma_1) &
 \Sigma_1^{-1}V_1^TB^TU_2 \\ U_2^TBV_1\Sigma_1^{-1} & A_{22}} U^H.$$ This shows that  $$\|A_o\|_F=\sqrt{\|D\circ (\Sigma_1 U_1^TBV_1 + V_1^TB^TU_1 \Sigma_1) \|_F^2+2\|U_2^TBV_1\Sigma_1^{-1}\|_F^2 + \|A_{22}\|_F^2}$$  is minimized when $A_{22}=0$ which yields the desired results for the Frobenius norm.

 For spectral norm, note that $\|A_o\|_2 \geq \mu_1,$ where
 $$\mu_1=\left\| \bmatrix{D\circ (\Sigma_1 U_1^TBV_1 + V_1^TB^TU_1 \Sigma_1) \\ U_2^TBV_1\Sigma_1^{-1} }\right\|_2.$$ By Theorem \ref{dkw}, we have $\|A_o\|_2=\mu_1$ when $$A_{22}=-K  ( \overline{ D\circ (\Sigma_1 U_1^TBV_1 + V_1^TB^TU_1 \Sigma_1)}) K^T + \mu_1 (I-KK^H)Z(I-\overline{K} K^T)^{1/2}$$ where $$K=U_2^TBX^\dagger U_1\left(\mu_1^2I -  (\overline{ D\circ (\Sigma_1 U_1^TBV_1 + V_1^TB^TU_1 \Sigma_1)})  (D\circ (\Sigma_1 U_1^TBV_1 + V_1^TB^TU_1 \Sigma_1))  \right)^{-1/2}$$ and $Z$ is an arbitrary contraction such that $Z=Z^T.$ Hence the desired form of $A_o$ follows for the spectral norm.

 The proof is similar when $\S=\herm$ and follows from Theorem~\ref{silsp:sym} and (\ref{eqn2:herm}). So, suppose that $\S=\sksym.$ Then by Theorem~\ref{silsp:sym} and (\ref{eqn2:sksym}), we have $$A_o =\overline{U}
 \bmatrix{D\circ (\Sigma_1 U_1^TBV_1 - V_1^TB^TU_1 \Sigma_1) &
 -\Sigma_1^{-1}V_1^TB^TU_2 \\ U_2^TBV_1\Sigma_1^{-1} & A_{22}} U^H.$$ As before, setting $A_{22} =0,$ the desired results follow for the Frobenius norm. On the other hand,  applying Theorem \ref{dkw} to the matrix $A_o$ and following steps similar to those in the case when $\S=\sym,$ we obtain the desired results for the spectral norm. $\blacksquare$




\begin{remark} We mention that the inverses in the expression of $K$ in Theorem~\ref{cor:silsp:sym} are replaced by their Moore-Penrose pseudo-inverses when the matrices are singular.
\end{remark}

%
%

Now, we consider the special case of a pair of nonzero vectors $x$ and $b$ in $ \K^n.$ It is well known~\cite{AlaA, fran6} that there always exists a symmetric matrix $A$ such that $Ax = b.$ Consequently, we only need to consider {\sc Problem-I} for skew-symmetric and  Hermitian matrices. We have the following result which follows from Theorems~\ref{silsp:sym} and~\ref{cor:silsp:sym}.

\begin{theorem}  Let $x$ and $b$ be nonzero vectors in $\K^n$ and  $\S \in \{\sksym, \herm\}.$ Let $P_x := I - \dfrac{xx^H}{\|x\|_2^2}.$  Then  we have
$$\rho^{\S}(x, b) =\left\{%
\begin{array}{ll}
    \dfrac{| x^Tb |}{\|x\|_2}, & \hbox{if\,  $\S=\sksym,$} \\ & \\
    \dfrac{|\im ( x^Hb) |}{\|x\|_2}, & \hbox{if \, $\S=\herm,$} \\
\end{array}%
\right.$$ Further, we have  $A \in \S$ and $\|Ax -b\|_F = \rho^{\S}(x, b)$ if and only if $A$ is of the form   $$A =\left\{%
\begin{array}{lll}
   \dfrac{1}{\|x\|_2^2} [P_x^Tb \bar x^T - \overline{x}b^TP_x] + P_x^T Z P_x, & \mbox{ if } \,\, \S=\sksym, \\ &  & \\
   \dfrac{\re(x^Hb)}{\|x\|^4_2} xx^H + \dfrac{1}{\|x\|^2_2}
[xb^HP_x + P_xbx^H] + P_x Z P_x,  &  \mbox{ if } \,\, \S=\herm, \\
\end{array}%
\right.$$ for some $Z\in \S.$

\von
{\bf (a) Frobenius norm: } We have $$\sigma^{\S}(x,b) = \left\{\begin{array}{ll} \sqrt{2}\|\bar P_x b\|_F = \sqrt{2} \sqrt{\|b\|_2^2 - |x^Tb|^2 /\|x\|_2^2 }, &  \mbox{ if } \S = \sksym, \\ & \\
\sqrt{2}\|P_x b\|_F = \sqrt{2} \sqrt{\|b\|_2^2 - |x^Hb|^2 /\|x\|_2^2 }, & \mbox{ if } \S = \herm.\end{array} \right.
$$
Further, there is a unique matrix $A_o \in \S$  such that  $\|Ax -b\|_F = \rho^{\S}(x, b)$ and $\|A_o\|_F=\sigma^{\S}(x,b),$ and is given by
$$A_o =  \left\{\begin{array}{ll} \dfrac{1}{\|x\|_2^2} [P_x^Tb \bar x^T - \overline{x}b^TP_x], & \mbox{ if } \S = \sksym, \\ & \\
\dfrac{\re(x^Hb)}{\|x\|^4_2}xx^H + \dfrac{1}{\|x\|_2^2} [xb^HP_x + P_xbx^H], & \mbox{ if } \S = \herm. \end{array} \right.
 $$

{\bf (b) Spectral norm:} We have
$$\sigma^{\S}(x,b) =  \left\{\begin{array}{ll}\|\bar P_x b\|_2 =\sqrt{\|b\|_2^2 - |x^Tb|^2 /\|x\|_2^2 } =:\mu, & \mbox{ if } \S = \sksym,\\ &\\
 \|P_x b\|_2 =\sqrt{\|b\|_2^2 - |x^Hb|^2 /\|x\|_2^2 } =:\mu, & \mbox{ if } \S = \herm. \end{array} \right.
$$
Further, we have  $A_o\in \S$  such that  $\|Ax -b\|_F = \rho^{\S}(x, b)$ and $\|A_o\|_2=\sigma^{\S}(x,b)$ if and only if $A_o$ is of the form
$$A_o =  \left\{\begin{array}{l} \dfrac{1}{\|x\|^2_2} [P_x^T b \bar x^T - \overline{x}b^TP_x] + \mu \overline{Q}_1(I-KK^H)^{1/2}Z(I-\overline{K}K^T)^{1/2}Q_1^H, \\ \hfill{ \mbox{ if } \S = \sksym, }\\ \\
 \dfrac{\re(x^Hb)}{\|x\|^4_2}xx^H + \dfrac{1}{\|x\|_2^2} [xb^HP_x + P_xbx^H] - \dfrac{\re(x^Hb) P_xbb^HP_x}{\mu^2 \|x\|_2^4 - (\re(x^Hb))^2} +\\ \\
  \hfill{ \mu Q_1 (I-KK^H)^{1/2}Z(I-KK^H)^{1/2}   Q_1^H,}  \hfill{ \mbox{ if } \S = \herm,} \end{array} \right.
$$ where  $Z\in\S$ is a contraction, $ Q_1\in\C^{n \times (n-1)}$ is an isometry such that $Q_1^Hx =0,$   $K = \dfrac{Q_1^Tb}{\mu\|x\|_2}$ when $ \S= \sksym,$  and   $K = \dfrac{Q_1^Hb}{\|x\|_2} \left(\mu^2 - \dfrac{\re(x^Hb)^2 }{\|x\|_2^4}\right)^{-1/2}$ when $\S= \herm.$ \end{theorem}

\vone{\bf Conclusion.} We have provided a complete solution of the structured inverse least-squared problem $\min_{A \in \S}\|AX-B\|_F$ when $\S$ is either a Jordan algebra or a Lie algebra associated with an appropriate scalar product. We have characterized (Theorem~\ref{silsp:sym}) solutions of the SILSP and have determined all optimal solutions (Theorem~\ref{cor:silsp:sym}) of the SILSP.

\end{document}